\documentclass[a4paper, reqno]{amsart}
\usepackage{amssymb, amsmath, amsthm, enumerate, mathabx, mathrsfs}

\numberwithin{equation}{section}
\newtheorem{thm}{Theorem}
\newtheorem{theo}[thm]{Theorem}

\title[A Variant of Yano's Extrapolation Theorem on Hardy Spaces]{A Variant of Yano's Extrapolation Theorem on Hardy Spaces}

\author[Odysseas Bakas]{Odysseas Bakas}
\address{Department of Mathematics, Stockholm University, 106 91 Stockholm, Sweden}

\email{bakas@math.su.se}



\date{}

\begin{document}

\maketitle

\begin{abstract} In this note  we prove a variant of Yano's classical extrapolation theorem for sublinear operators acting on analytic Hardy spaces over the torus. 
\end{abstract}

\section{Introduction}

Let $(X, \mu) $ and $(Y, \nu)$ be two finite measure spaces.  If $T : (X, \mu) \rightarrow (Y, \nu ) $ is a sublinear operator such that there exist constants $C_0, r >0$ satisfying 
\begin{equation}\label{op_norm}
\sup_{\| g \|_{L^p (X)} =1 } \| T (g)  \|_{ L^p (Y) } \leq C_0 (p-1)^{-r} ,
 \end{equation}
 for every $ 1 < p \leq 2 $, then a classical theorem of S. Yano \cite{Yano} asserts that 
\begin{equation}\label{Orlicz_bound}
 \| T (f)  \|_{ L^1 (Y) } \leq A + B \int_{X} |f(x)| \log^r (1+ |f(x)|) d\mu(x)
 \end{equation}
 for all simple functions $f$ on $X$,  where $A,B >0$ are constants depending only on $C_0$, $r$, $ \mu(X)$, and $\nu (Y)$. See also Theorem 4.41 in Chapter XII of \cite{Zygmund_book}. 

In this note we prove a version of the aforementioned extrapolation theorem of Yano for sublinear operators acting on functions belonging to analytic Hardy spaces over the torus.  Namely, we prove that if $T$ is a sublinear operator acting on functions defined over the torus such that its  operator norm from $(H^p (\mathbb{T}), \| \cdot \|_{L^p (\mathbb{T})})$ to $(L^p (\mathbb{T}), \| \cdot \|_{L^p (\mathbb{T})})$ behaves like $(p-1)^{-r}$ as $p \rightarrow 1^+$, then $T$ satisfies an inequality analogous to (\ref{Orlicz_bound}) for functions in  Hardy spaces, see Theorem \ref{main_theorem} below. Here, for $1\leq p \leq \infty$, $H^p (\mathbb{T}) $ denotes the analytic Hardy space $H^p$ on $\mathbb{T}$ given by  
$$ H^p (\mathbb{T}) = \{ f \in L^p (\mathbb{T}) : \widehat{f} (n) = 0 \ \mathrm{for} \ n <0\}  .$$

The  study of such variants of Yano's extrapolation theorem in the present paper is motivated by some classical results of S. Pichorides \cite{Pichorides} and A. Zygmund \cite[Theorem 8]{Zygmund} on mapping properties of the Littlewood-Paley square function $S$ ``near'' $H^1 (\mathbb{T})$. Recall that given a trigonometric polynomial $f$ on $\mathbb{T}$, the classical Littlewood-Paley square function $S(f)$ of $f$ is defined as 
$$ S (f) := \Big( \sum_{k \in \mathbb{Z} } |\Delta_k (f) |^2 \Big)^{1/2}, $$ 
where   $\Delta_0 (f) (\theta) := \widehat{f} (0)$,  $\theta \in \mathbb{T}$ and for $k \in \mathbb{N}$,
$$ \Delta_k (f) (\theta) := \sum_{n=2^{k-1}}^{2^k - 1} \widehat{f} (n) e^{ i 2 \pi n \theta} \ \mathrm{and}\  \Delta_{-k} (f) (\theta) := \sum_{n=-2^{k} + 1}^{-2^{k - 1}} \widehat{f} (n) e^{ i 2 \pi n \theta} $$ 
for $\theta \in \mathbb{T}$. In \cite{Zygmund}, Zygmund showed that for every function $f \in H^1 (\mathbb{T})$ one has
\begin{equation}\label{Zyg_ineq}
\| S (f) \|_{L^1 (\mathbb{T})} \lesssim 1 +  \int_{\mathbb{T}} |f(\theta)| \log (1+ |f(\theta)|) d\theta   \ \ (f \in H^1 (\mathbb{T}))
\end{equation}
and in \cite{Pichorides}, Pichorides showed that
\begin{equation}\label{Pich_thm}
\sup_{\substack{  g \in H^p (\mathbb{T}) : \\ \| g \|_{L^p (\mathbb{T})} = 1 } }\| S (g) \|_{L^p (\mathbb{T})} \sim (p-1)^{-1}
\end{equation}
as $p \rightarrow 1^+$. Both of these results were originally proved by using one-dimensional complex-analytic techniques such as the analytic factorisation of Hardy spaces on the torus. 
Recently, in \cite{BRS}  the aforementioned results of Pichorides and Zygmund were extended to higher dimensions in a  unified way using ``real-variable'' techniques. In particular, it is shown in \cite{BRS} that both (\ref{Zyg_ineq}) and (\ref{Pich_thm}) can be obtained using a result of T. Tao and J. Wright on endpoint mapping properties of Marcinkiewicz multiplier operators \cite{TW} together with a ``Marcinkiewicz-type'' interpolation argument for $H^p (\mathbb{T})  $ spaces due to S. Kislyakov and Q. Xu \cite{KX_2}, see also \cite{Bourgain_1}.
Therefore, motivated by the above remarks, one is naturally led to ask whether it is possible to deduce the result of Zygmund (\ref{Zyg_ineq}) directly from that of Pichorides (\ref{Pich_thm})  by using some ``real-variable'' approach or, even more generally, whether a general extrapolation theorem for sublinear operators acting on $H^p (\mathbb{T})$ spaces holds true. In this note we prove that this is indeed the case, namely we have the following result.

\begin{theo}\label{main_theorem}
Let $T$ be a sublinear operator acting on functions defined over $\mathbb{T}$, namely for all measurable functions $f, g$ on $\mathbb{T}$ and each complex number $\alpha$ one has $|T (f+g) | \leq |T (f)| + |T (g)|$ and $|T (\alpha f) | = | \alpha | |T(f)|$. 

If there exist constants $C_0, r>0$ such that
\begin{equation}\label{norm_Hardy}
 \sup_{\substack{  g \in H^p (\mathbb{T}) : \\ \| g \|_{L^p (\mathbb{T})} = 1 } }\| T (g) \|_{L^p (\mathbb{T})} \leq C_0 (p-1)^{-r} 
 \end{equation}
 for every $1 < p \leq 2$, then there exists a constant $D>0$, depending only on $C_0, r$, such that
\begin{equation}\label{Orlicz_Hardy}
 \| T (f) \|_{L^1 (\mathbb{T})} \leq D  \| f \|_{L \log^r L (\mathbb{T}) } 
 \end{equation}
for every analytic trigonometric polynomial $f$ on $\mathbb{T}$.
\end{theo}

Note that if $T$ satisfies the assumptions of Theorem \ref{main_theorem}, then  it is easy to see that for every analytic trigonometric polynomial $f$ one has $  \| T (f) \|_{L^1 (\mathbb{T})} \lesssim  \| f \|_{L \log^{r+1} L (\mathbb{T}) }  $. Indeed, if $P$ denotes the Riesz projection onto non-negative frequencies, namely $P$ is the multiplier operator on $\mathbb{T}$ with symbol $\chi_{\mathbb{N}_0}$, then it is a standard fact that $\| P \|_{L^p (\mathbb{T}) \rightarrow L^p (\mathbb{T})} \lesssim (p-1)^{-1}$ as $p \rightarrow 1^+$. Hence, if $T$ is a sublinear operator satisfying (\ref{norm_Hardy}), then $ \| T \circ P \|_{L^p (\mathbb{T}) \rightarrow L^p (\mathbb{T})} \lesssim (p-1)^{-(r+1)}$ as $p \rightarrow 1^+$. It thus follows from Yano's extrapolation theorem applied to $T \circ P$ that
$T$ satisfies the aforementioned weaker version of (\ref{Orlicz_Hardy}) where $\| f \|_{L \log^r L (\mathbb{T})}$ is replaced by $\| f \|_{L \log^{r+1} L (\mathbb{T})}$, $f$ being an analytic trigonometric polynomial. 
Therefore, for sublinear operators satisfying (\ref{norm_Hardy}), Theorem \ref{main_theorem} improves the trivial exponent $s=r+1$ in $L \log^s L (\mathbb{T})$ to the optimal one; $s =r$, see discussion in Subsection \ref{optimality}. 

The proof of our result is based on Yano's original argument  \cite{Yano}  
combined with some well-known techniques on interpolation between $H^p (\mathbb{T})$ spaces, see  \cite{Bourgain_2}, \cite{K}.
 At this point, it is worth mentioning that the main idea in Yano's paper is to decompose a given function $f$ as
$$ f = \sum_{n \in \mathbb{N}_0 } f_n,$$
where $ f_0 = \chi_{\{ |f| < 1 \} } f$ and $f_n = \chi_{\{ 2^{n-1} \leq |f| < 2^n \} } f$ (for $n \in  \mathbb{N}$) and then apply the assumption (\ref{op_norm}) to each $f_n$ separately for $p = p_n = 1 +1/(n+1)$, noting that $f_n \in L^{p_n} (\mathbb{T})$ for $n \in \mathbb{N}_0 $. Unfortunately, given a function $f \in H^1 (\mathbb{T})$, if one defines $f_n$ as above, then $f_n$ are not necessarily in the analytic Hardy space $H^{p_n} (\mathbb{T})$ anymore and hence, one cannot apply (\ref{norm_Hardy}) to $f_n$.  To surpass this difficulty, the idea is to use an ``analytic decomposition of unity'' of $f = \sum_n \widetilde{f}_n$ which is due to Kislyakov \cite{K}, see also J. Bourgain's paper \cite{Bourgain_2}. In particular, each $\widetilde{f}_n$ may be regarded as an appropriate ``$H^{\infty}$-replacement'' of $f_n$ in the sense that $\widetilde{f}_n \in H^{\infty} (\mathbb{T})$, $f = \sum_n \widetilde{f}_n$ and each $\widetilde{f}_n$ ``essentially behaves like'' $f_n$, $f_n $ being as above. 
The proof of Theorem \ref{main_theorem} is given in Section \ref{main_section}. In Section \ref{discussion} we briefly discuss about some further remarks related to the present work.

\subsection*{Notation} 
We denote the set of integers by $\mathbb{Z}$. The set of natural numbers is denoted by $\mathbb{N}$ and the set of non-negative integers is denoted by $\mathbb{N}_0$. 
 
We identify functions over the torus $\mathbb{T}$ with functions defined on the set $[0,1)$. 

If $f \in L^1 (\mathbb{T})$ is such that $\mathrm{supp} (  \widehat{f}) $ is finite, then $f$ is said to be a trigonometric polynomial on $\mathbb{T}$. If $f$ is a trigonometric polynomial on $\mathbb{T}$ such that $\mathrm{supp}(\widehat{f}) \subset \mathbb{N}_0$, then we say that $f$ is an analytic trigonometric polynomial on $\mathbb{T}$.

For $r>0$, $L \log^r L (\mathbb{T})$ denotes the class of all measurable functions $f$ on $\mathbb{T}$ satisfying $\int_{\mathbb{T}} |f (\theta)| \log^r (1 + |f (\theta) |) d \theta < \infty $. For $f \in L \log^r L (\mathbb{T})$, if we set 
$$ \| f \|_{L \log^r L (\mathbb{T})} := \inf \Big\{  \lambda >0 : \int_{\mathbb{T}} \Phi_r (\lambda^{-1} |f (\theta)| ) d \theta \leq 1 \Big\},  $$
where $\Phi_r (x) = x ( [1+\log (x+1)]^r - 1 )$ $(x \geq 0)$, then $ \| \cdot \|_{L \log^r L (\mathbb{T})} $ is a norm on $L \log^r L (\mathbb{T})$ and, moreover,  $(L \log^r L (\mathbb{T}), \| \cdot \|_{L \log^r L (\mathbb{T})} )$ is a Banach space. For more details on Orlicz spaces, see \cite{KR}.

If $x$ is a real number, then $\lfloor x \rfloor $ denotes its integer part. As usual, $\log x$ denotes the natural logarithm of a positive real number $x$. The logarithm of $x>0$ to the base $2$ is denoted by $\log_2 (x)$.

Given two positive quantities $X$ and $Y$, if there exists a positive constant $C>0$ such that $X \leq C Y$ we shall write $X \lesssim Y$. If the constant $C$ depends on some parameters $s_1, \cdots, s_n$ then we shall also write $ X \lesssim_{s_1, \cdots, s_n} Y$. Moreover, if $X \lesssim Y$ and $Y \lesssim X$ we write $X \sim Y$.

\section{Proof of Theorem \ref{main_theorem}}\label{main_section} Let $T$ be a sublinear operator satisfying (\ref{norm_Hardy}) and let $f$ be a fixed analytic trigonometric polynomial on $ \mathbb{T}$. We shall prove that
\begin{equation}\label{Orlicz_analytic}
 \| T (f) \|_{L^1 (\mathbb{T})} \leq A + B  \int_{\mathbb{T}} |f (\theta)| \log^r (1 + |f (\theta)|)  d \theta  ,
\end{equation}
where $A,B>0$  depend only on $C_0, r$ and not on $f$. Towards this aim, following \cite{KX_1} (see also \cite{K}), for $\lambda >0$ consider the function
$$ a_{\lambda} (\theta) := \max \Big\{ 1, \Big( \frac{|f(\theta)|}{\lambda} \Big)^{1/3} \Big\} $$
and then define
$$ F_{\lambda} ( \theta) := \frac{1}{a_{\lambda} (\theta) + i H (a_{\lambda} ) (\theta)}$$
and
$$ G_{\lambda} (\theta):= 1 - ( 1 - [F_{\lambda} (\theta)]^4 )^4. $$
Here, $H$ denotes the periodic Hilbert transform. It follows that $F_{\lambda} , G_{\lambda} \in H^{\infty} (\mathbb{T})$, see e.g. the proof of  \cite[Lemma 7.4.2]{Pavlovic}. Moreover, since $| F_{\lambda}| \leq \min\{ 1, \lambda^{1/3} |f|^{-1/3} \} $, one has
\begin{equation}\label{characteristic}
 |G_{\lambda}| = | 1 - (1 - [ F_{\lambda} ]^4 )^4 | \leq A_0  | F_{\lambda} |^4   \leq  A_0 | F_{\lambda} |^3 \leq A_0 \min\{ 1, \lambda |f|^{-1 } \} ,
 \end{equation}
 where $A_0 >0$ is an absolute constant. In particular, $|G_{\lambda} f   | \lesssim \lambda$ on $\mathbb{T}$ and so, in order to define an appropriate ``bounded analytic replacement'' of $\chi_{\{ \lambda \leq |f| < 2 \lambda \} } f$, one is led to consider functions of the form $ (G_{2 \lambda} - G_{\lambda}) f$.  More precisely, arguing as in the proof of  \cite[Lemma 4.2]{KX_1}, consider the functions $(\widetilde{f}_n)_{n \in \mathbb{N}_0}$ in $H^{\infty} (\mathbb{T})$ given by
$$ \widetilde{f}_0 : = G_1 f$$
and
$$ \widetilde{f}_n := (G_{2^n} - G_{2^{n-1}} ) f $$
for $n \in \mathbb{N}$. Note that  there exists an $N \in \mathbb{N}$, depending on $f$, such that $\widetilde{f}_k \equiv 0 $ for all $k \geq N$. Indeed, since $f$ is an analytic trigonometric polynomial, if we take $N \in \mathbb{N}$ such that $2^{N-1} > \| f \|_{L^{\infty} (\mathbb{T})}$, then for every $k \geq N $ one has $|f (\theta)| \leq 2^{k-1}$ for all $\theta \in \mathbb{T} $. Hence, for every $k \geq N$ one has $a_{2^{k-1}}   = a_{2^k} \equiv 1$ on $ \mathbb{T} $. Therefore,  $G_{2^{k-1}}   = G_{2^k} \equiv 1$ on $  \mathbb{T} $ and we thus deduce that $\widetilde{f}_k = (G_{2^k} - G_{2^{k-1}}) f \equiv 0$ on $\mathbb{T}$ whenever $k \geq N$. So, one has the decomposition  
 $$ f (\theta) = \sum_{n=0}^N \widetilde{f}_n  (\theta)$$
 for all $\theta \in \mathbb{T} $.  Next, as in Yano's paper \cite{Yano}, using the sublinearity of $T$ and then H\"older's inequality, one  deduces that
$$   \| T (f) \|_{L^1 (\mathbb{T} )} \leq \sum_{n=0}^N \| T (\widetilde{f}_n) \|_{L^{p_n} (\mathbb{T} )}, $$
where $p_n = 1 + 1/(n+1)$. Hence, using our assumption (\ref{norm_Hardy}) for each $ 0 \leq n \leq N$, one gets
\begin{equation}\label{step_1}
 \| T (f) \|_{L^1 (\mathbb{T})} \leq  C_0 \sum_{n=0}^N (n+1)^r \| \widetilde{f}_n  \|_{L^{p_n} (\mathbb{T})}.
 \end{equation}
Note that since
$$ | \widetilde{f}_n | \leq |G_{2^{n}}| |f| + |G_{2^{n-1}} | |f| $$
it follows from (\ref{characteristic}) that
\begin{equation}\label{est_1}
|\widetilde{f}_n| \leq A'_0 2^n,
\end{equation}
 where $A'_0 >0 $ is an absolute constant, independent of $f$ and $n$.  Hence, we have 
 \begin{align*}
 \| \widetilde{f}_n  \|_{L^{p_n} (\mathbb{T})} &= A'_0 2^n \Big( \int_{ \mathbb{T} } \big[  (A'_0 2^n)^{-1} | \widetilde{f}_n (\theta)| \big]^{\frac{n+2}{n+1} } d \theta \Big)^{ \frac{n+1}{n+2} } \\
 & \leq A'_0 2^n   \Big( \int_{ \mathbb{T} }   (A'_0 2^n)^{-1} | \widetilde{f}_n (\theta)|  d \theta \Big)^{ \frac{n+1}{n+2} } 
 \end{align*}
 and we thus deduce that
 $$ \| \widetilde{f}_n  \|_{L^{p_n} ( \mathbb{T} )} \lesssim \Big( \int_{ \mathbb{T} }    | \widetilde{f}_n (\theta)|  d \theta \Big)^{ \frac{n+1}{n+2} } , $$
 where the implied constant is independent of $n$. Using now the elementary inequality 
 $$ t^{(n+1)/(n+2)} \leq e^{r+2} t + (n+1)^{-(r+2)}$$
which is valid for all $t > 0$ and $n \geq 0$, we obtain
\begin{equation}\label{est_2}
 \| \widetilde{f}_n  \|_{L^{p_n} ( \mathbb{T} )} \lesssim_r   \int_{\mathbb{T}}   | \widetilde{f}_n (\theta)| d \theta  + \frac{1}{(n+1)^{r+2}},
 \end{equation}
where the implied constant depends only on $r$ and not on $f, n$.  Moreover, note that by using (\ref{est_1}) for $n=0$, one has
\begin{equation}\label{trivial}
\| \widetilde{f}_0 \|_{L^1 ( \mathbb{T} )} \lesssim 1 .
\end{equation}
Hence, it follows from (\ref{step_1}), (\ref{est_2}), and (\ref{trivial}) that
$$
\| T (f ) \|_{L^1 ( \mathbb{T} )} \lesssim_{C_0, r} 1 + \sum_{n=1}^N (n+1)^r  \| \widetilde{f}_n \|_{L^1 ( \mathbb{T} )}. $$
Therefore, to prove (\ref{Orlicz_analytic}) it suffices to show that
\begin{equation}\label{main_estimate}
 \sum_{n=0}^N (n+1)^r  \| \widetilde{f}_n \|_{L^1 ( \mathbb{T} )}  \lesssim_r 1+\int_{ \mathbb{T} } |f (\theta)| \log^r (1 + |f (\theta)|) d \theta.
\end{equation}
To this end, the idea is to write
\begin{align*}
 \| \widetilde{f}_n \|_{L^1 ( \mathbb{T} )} &=  \int_{\{ |f| < 2^{n-1} \}} |\widetilde{f}_n (\theta)| d \theta +  \int_{\{ |f| \geq 2^{n-1} \}} |\widetilde{f}_n (\theta)| d \theta \\
 & =  I_n^{(1)} + I_n^{(2)}  
 \end{align*}
and then treat each of the terms $I_n^{(1)}$, $I_n^{(2)}$  separately showing that
\begin{equation}\label{sufficient}
  \sum_{n=1}^N (n+1)^r I_n^{(i)} \lesssim_r 1 + \int_{ \mathbb{T} } |f (\theta)| \log^r (1 + |f (\theta)|) d \theta \ \ (i =1,2).  
\end{equation}
Note that (\ref{sufficient}) (for $i=1,2$) automatically establishes (\ref{main_estimate}). We have thus reduced matters to showing  (\ref{sufficient}) for $i=1,2$. 

The case $i=2$ in (\ref{sufficient}) is the easiest one; using (\ref{est_1}) for each   $1 \leq n \leq N+1$, we get
\begin{align*}
\sum_{n=1}^N (n+1)^r I_n^{(2)} &=  \sum_{n=1}^N (n+1)^r \int_{\{ |f| \geq 2^{n-1} \}} |\widetilde{f}_n (\theta)| d \theta \\
&\lesssim \sum_{n=1}^N (n+1)^r \int_{\{ |f| \geq 2^{n-1} \}} 2^n  d \theta .
\end{align*} 
Hence, by using Fubini's theorem, we obtain
\begin{align*}
\sum_{n=1}^N (n+1)^r I_n^{(2)} &\lesssim \int_{ \mathbb{T} } \Big( \sum_{n=1}^{2 + \lfloor \log_2 (|f (\theta)|+1) \rfloor  } 2^n (n+1)^r  \Big) d\theta \\ 
&\lesssim_r 1+ \int_{ \mathbb{T} } |f(\theta)| \log^r (1+|f (\theta)|) d \theta,
 \end{align*} 
as desired.  

It remains to prove (\ref{sufficient}) for $i=1$. Towards this aim, note that since
$$ G_{2^n} - G_{2^{n-1}}  = (1 - [F_{2^n-1}]^4)^4 -  (1 - [F_{2^n}]^4)^4$$
and $|F_{2^n-1}|  \leq 1$, $ |F_{2^n}|  \leq 1 $ on $\mathbb{T}$, it follows that
$$ | G_{2^n} - G_{2^{n-1}} | \leq |1 - [F_{2^n}]^4|^4 + |1 - [F_{2^n-1}]^4|^4  \lesssim |1 - F_{2^n}|^2 +  |1 - F_{2^{n-1}}|^2 $$
on $\mathbb{T}$, where the implied constant is independent of $n$. Hence, we have
\begin{align*}
 I^{(1)}_n &=   \int_{\{ |f| < 2^{n-1}\}} |\widetilde{f}_n (\theta)| d \theta \\
 & \lesssim \int_{\{ |f| < 2^{n-1}\}}  |1 - F_{2^n} (\theta) |^2 |f (\theta) | d \theta + \int_{\{ |f| < 2^{n-1}\}}  |1 - F_{2^{n-1}} (\theta) |^2 |f (\theta) |  d \theta \\
 & \leq  I^{(1, \alpha)}_n + I^{(1, \beta)}_n,
 \end{align*}
where
$$ I^{(1, \alpha)}_n  := 2^{n-1}\int_{\{ |f| < 2^{n-1}\}}  |1 - F_{2^n} (\theta) |^2  d \theta$$
and
$$ I^{(1, \beta)}_n  := 2^{n-1}\int_{\{ |f| < 2^{n-1}\}}  |1 - F_{2^{n-1}} (\theta) |^2  d \theta .$$
We shall prove that
\begin{equation}\label{eq_1}
 I^{(1, \alpha)}_n \lesssim 2^{n/3} \int_{|\{ f| \geq 2^n \} } |f (\theta)|^{2/3} d \theta
\end{equation}
and
\begin{equation}\label{eq_2}
 I^{(1, \beta)}_n \lesssim 2^{n/3} \int_{ \{ |f| \geq 2^n \} } |f (\theta)|^{2/3} d \theta.
\end{equation}

To prove (\ref{eq_1}), we shall again make use of arguments from \cite{KX_1} (see also \cite{K}); since the periodic Hilbert transform $H$ of any constant function on $ \mathbb{T} $ is identically $0$, one may write
$$ I^{(1, \alpha)}_n = 2^{n-1}   \int_{\{ |f| < 2^{n-1}\}}  \Big| \frac{(1-a_{2^n} (\theta)) + i H (1- a_{2^n})(\theta)}{a_{2^n} (\theta) + i H (1- a_{2^n})(\theta)} \Big|^2 d \theta $$
and since $a_{2^n} \equiv 1$ on $\{ |f| < 2^{n-1}\}$, it follows that
$$   I^{(1, \alpha)}_n  \leq 2^{n-1}  \int_{\{ |f| < 2^{n-1}\}}  | H  (1- a_{2^n})(\theta) |^2 d \theta \leq  2^{n-1}  \| H  (1- a_{2^n}) \|^2_{L^2 ( \mathbb{T} )}. $$
Hence, by using the $L^2$-boundedness of $H $, one gets
\begin{align*}
 I^{(1, \alpha)}_n \leq 2^{n-1}  \| 1- a_{2^n}  \|^2_{L^2 ( \mathbb{T} )} &= 2^{n-1}  \int_{\{ | f  | \geq 2^n \} } \Big| \Big( \frac{|f (\theta)|}{2^n} \Big)^{1/3} - 1\Big|^2 d \theta \\
& \leq 2^{n/3} \int_{\{ | f | \geq 2^n\} } |f (\theta)|^{2/3} d \theta
\end{align*}
and this completes the proof of (\ref{eq_1}). The proof of (\ref{eq_2}) is completely analogous. It thus follows from (\ref{eq_1}) and (\ref{eq_2}) that
$$ I^{(1)}_n \lesssim 2^{n/3} \int_{\{ | f | \geq 2^n\} } |f (\theta)|^{2/3} d \theta $$
and hence, by using Fubini's theorem, we have
\begin{align*}
\sum_{n=1}^N (n+1)^r I_n^{(1)} &\lesssim \int_{ \mathbb{T} } | f (\theta)|^{2/3} \Big( \sum_{n=1}^{1 +  \lfloor \log_2 (|f (\theta)|+1) \rfloor } 2^{n/3} (n+1)^r  \Big) d\theta \\ 
&\lesssim_r 1+ \int_{ \mathbb{T} } |f(\theta)| \log^r (1+|f (\theta)|) d \theta.
 \end{align*} 
Therefore, the proof of (\ref{sufficient}) for $i=1$ is complete. So, (\ref{sufficient}) holds for $i=1,2$ and hence, (\ref{main_estimate}) is also true. We have thus shown that (\ref{Orlicz_analytic}) holds for every analytic trigonometric polynomial on the torus.

 To complete the proof of Theorem \ref{main_theorem}, note that  (\ref{Orlicz_Hardy}) follows directly from  (\ref{Orlicz_analytic}) using a simple scaling argument. Indeed, if $f$ is an analytic trigonometric polynomial on $\mathbb{T}$ with $\| f \|_{L \log^r L (\mathbb{T})} =1 $, then $\int_{\mathbb{T} } | f (\theta)| \log^r (1+|f (\theta)|) d \theta \lesssim_r 1$. Hence, one deduces from (\ref{Orlicz_analytic}) that if $\| f \|_{L \log^r L (\mathbb{T})} =1$ then $ \| T (f) \|_{L^1 (\mathbb{T})} \leq D $, where $D>0$ is an absolute constant depending only on $C_0, r$. In the general case, if $f$ is a non-zero analytic trigonometric polynomial on $\mathbb{T}$ then, using the previous implication and the scaling invariance of $T$,  (\ref{Orlicz_Hardy}) follows.

\section{Some Further Remarks}\label{discussion} 

\subsection{Sharpness of Theorem \ref{main_theorem}}\label{optimality}
As mentioned in the introduction, Zygmund's inequality (\ref{Zyg_ineq}) on the classical Littlewood-Paley square function can now be obtained as a corollary of Pichorides's result (\ref{Pich_thm}) via Theorem \ref{main_theorem}. Moreover, the exponent $r = 1$ in the Orlicz space $L \log^r L (\mathbb{T})$ in (\ref{Zyg_ineq})  cannot be improved, see \cite{BRS}. Therefore, the example of the Littlewood-Paley square function $S$ shows that Theorem \ref{main_theorem} is in general sharp in the following sense; there exists a sublinear operator $T$ and an $r = r (T) >0$ such that $\| T \|_{(H^p (\mathbb{T}), \| \cdot \|_{L^p(\mathbb{T})})  \rightarrow (L^p (\mathbb{T}), \| \cdot \|_{L^p(\mathbb{T})}) } \sim (p-1)^{-r}$ as $p \rightarrow 1^+$ and $ \| T (f) \|_{L^1 (\mathbb{T})} \lesssim \| f \|_{L \log^s L (\mathbb{T}) }$
holds for every analytic trigonometric polynomial $f$  when $ s = r$ but it \emph{does not} hold for any exponent $ s<r$. In particular, for the Littlewood-Paley square function $S$ one has $r = r(S) =1$. 

Note that if one removes the analyticity assumptions, then the behaviour of $S$ ``near'' $L^1 (\mathbb{T})$ is different than the one mentioned above. In particular,  Bourgain showed in \cite{Bourgain_3} that the $L^p (\mathbb{T}) \rightarrow L^p (\mathbb{T})$ operator norm of $S$ behaves like $(p-1)^{-3/2}$ as $p \rightarrow 1^+$. Moreover, $f \in L \log^{3/2} L (\mathbb{T})$ implies that $S (f) \in L^1 (\mathbb{T})$ and the exponent $r=3/2$ in $L \log^{3/2} L (\mathbb{T})$ is best possible, see \cite{Bakas}.

Another example illustrating Theorem \ref{main_theorem} and its sharpness (in the sense discussed above) is given by some classical results of Y. Meyer \cite{Meyer} and A. Bonami \cite{Bonami}. More specifically, consider the following subset of positive integers
$$\Lambda = \{ 3^k - 3^m : \ 0 \leq m \leq k-1 \ \mathrm{and}\ k,m \in \mathbb{Z} \} .$$
Then \cite[Corollaire 4]{Bonami} asserts that there exists an absolute constant $C >0$ such that  for every $q >2$ one has
\begin{equation}\label{Bonami_thm}
\| h \|_{L^q (\mathbb{T})} \leq C  q \| h \|_{L^2 (\mathbb{T})}
\end{equation}
 for all trigonometric polynomials $h$ on $\mathbb{T}$ with $\mathrm{supp}(\widehat{h}) \subset \Lambda$. Hence, if we consider the multiplier operator $T_{\Lambda} $  on $\mathbb{T}$ with symbol $\chi_{\Lambda}$, namely for every  trigonometric polynomial $f$  one has
$$ T_{\Lambda} (f) (\theta) = \sum_{n \in \Lambda} \widehat{f} (n) e^{i 2 \pi n \theta}$$
for all $\theta \in \mathbb{T}$,  then it follows from (\ref{Bonami_thm}) and duality that $T_{\Lambda}$ satisfies
\begin{equation}\label{Bonami_1}
\sup_{\substack{  g \in L^p (\mathbb{T}) : \\ \| g \|_{L^p (\mathbb{T})} = 1 } }\| T_{\Lambda} (g) \|_{L^p (\mathbb{T})} \lesssim (p-1)^{-1} \ \ (\mathrm{as}\ p \rightarrow 1^+)
\end{equation} 
and
\begin{equation}\label{Bonami_2}
\| T_{\Lambda} (f) \|_{L^1 (\mathbb{T}) } \lesssim \| f \|_{L \log L (\mathbb{T})}.
\end{equation} 
However, if we restrict ourselves to analytic Hardy spaces on the torus, then it follows from parts (a) and (c) of Th\'eor\`eme 1 on p. 549-550 in \cite{Meyer} that one has the improved bounds
\begin{equation}\label{Meyer_1}
\sup_{\substack{  g \in H^p (\mathbb{T}) : \\ \| g \|_{L^p (\mathbb{T})} = 1 } }\| T_{\Lambda} (g) \|_{L^p (\mathbb{T})} \lesssim (p-1)^{-1/2} \ \ (\mathrm{as}\ p \rightarrow 1^+)
\end{equation}
  and
 \begin{equation}\label{Meyer_2}
  \| T_{\Lambda} (f) \|_{L^1 (\mathbb{T})} \lesssim  \| f \|_{L \log^{1/2} L (\mathbb{T}) }  \ \ (f \in H^1 (\mathbb{T})),
 \end{equation}
 respectively. In particular, the linear operator $T_{\Lambda}$ satisfies (\ref{norm_Hardy}) and  (\ref{Orlicz_Hardy}) in Theorem \ref{main_theorem} for $r=1/2$.
 
Furthermore, the exponents $r=1/2$ in $(p-1)^{-1/2}$ in (\ref{Meyer_1}) and $r=1/2$ in $L \log^{1/2} L (\mathbb{T})$ in (\ref{Meyer_2}) cannot be improved. Indeed, to see that (\ref{Meyer_2}) is sharp, note that, by
using  (\ref{Bonami_thm}) and \cite[(1.4.1)]{Rudin}, one has
\begin{equation}\label{equivalence}
\| T_{\Lambda} (f) \|_{L^1 (\mathbb{T})} \leq \| T_{\Lambda} (f) \|_{L^2 (\mathbb{T})} \lesssim  \| T_{\Lambda} (f) \|_{L^1 (\mathbb{T})}
\end{equation} 
for every trigonometric polynomial $f $ on $\mathbb{T}$. Assume now that for some $r>0$ one has
\begin{equation}\label{test_r}
 \| T_{\Lambda} (f) \|_{L^1 (\mathbb{T})} \lesssim  \| f \|_{L \log^r L (\mathbb{T}) }  \ \ (f \in H^1 (\mathbb{T})). 
 \end{equation}
To show that $r \geq 1/2$, for a large $N \in \mathbb{N}$ that will eventually be sent to infinity, consider the function  $\beta_N \in H^{\infty} (\mathbb{T})$ given by
$$ \beta_N (\theta) = e^{i 2 \pi (2 \cdot 3^N +1) \theta} V_{3^N} (\theta), $$
where $V_n = 2 K_{2n+1} -  K_n$ is the de la Vall\'ee Poussin kernel of order $n \in \mathbb{N}$ and  $K_n$ denotes the $n$-th F\'ejer kernel, i.e. $K_n (\theta) = \sum_{|j| \leq n} [1-|j|/(n+1)] e^{i 2 \pi  j \theta}$, $\theta \in \mathbb{T}$. Since $ \| \beta_N \|_{L \log^r L (\mathbb{T})} \lesssim N^r $ and $\widehat{\beta_N} (j) =1 $ for all $j \in \mathbb{N}$ with $3^N  \leq j \leq  3^{N+1}+2 $, it follows from (\ref{test_r}) and  (\ref{equivalence}) that
\begin{align*}
N^r \gtrsim \| \beta_N \|_{L \log^r L (\mathbb{T})} \gtrsim \| T_{\Lambda} (\beta_N) \|_{L^1 (\mathbb{T})} &\gtrsim  \| T_{\Lambda} (\beta_N) \|_{L^2 (\mathbb{T})} \\
&= \Big( \sum_{k=1}^{N+2} \sum_{m=0}^{k-1}  | \widehat{\beta_N} (3^k - 3^m) |^2 \Big)^{1/2} \\
&\geq \Big(   \sum_{m=0}^N  | \widehat{\beta_N} (3^{N+1} - 3^m) |^2 \Big)^{1/2}  \sim N^{1/2}
\end{align*}
and so, by taking $N \rightarrow \infty$, one deduces that $ r \geq 1/2$, i.e. (\ref{Meyer_2}) is sharp. Similarly, one shows that the exponent $r=1/2$ in $(p-1)^{-1/2}$ in (\ref{Meyer_1}) as well as the exponents $r=1$ in $(p-1)^{-1}$ in (\ref{Bonami_1}) and $r = 1$ in $L \log L (\mathbb{T})$ in (\ref{Bonami_2}) are best possible.
\subsection{Extension of Theorem \ref{main_theorem} to Hardy-Orlicz spaces}\label{extension}
If $T$ satisfies the assumptions of Theorem \ref{main_theorem} and, moreover, $|T(f) -T (g)| \leq |T(f-g)|$ holds, then one can easily show that  (\ref{Orlicz_Hardy})  can be extended to all functions $f$ such that $\| f \|_{L \log^r L (\mathbb{T})} < \infty $ and $\mathrm{supp}(\widehat{f}) \subset \mathbb{N}_0$. To this end, one adapts e.g. the density argument on p. 120 in Chapter XII of \cite{Zygmund_book}. More precisely, fix a function $f $ with $\| f \|_{L \log^r L (\mathbb{T})} < \infty $ and $\mathrm{supp}(\widehat{f}) \subset \mathbb{N}_0$. By \cite[Proposition 3.4]{LQR}, there exists a sequence of analytic trigonometric polynomials $(\phi_n)_{n \in \mathbb{N}}$ that converges to $f$ in  $( L \log^r L (\mathbb{T}), \| \cdot \|_{L \log^r L (\mathbb{T})} )$. It thus follows from  (\ref{Orlicz_Hardy}) that the sequence $(T (\phi_n))_{n \in \mathbb{N}}$ is Cauchy in $(L^1 (\mathbb{T}) , \| \cdot \|_{L^1 (\mathbb{T})})$ and so, it converges to some $g \in L^1 (\mathbb{T})$ and $g$ is independent of the choice of $(\phi_n)_{n \in \mathbb{N}}$. Hence, if we set $T (f) := g$, we deduce that (\ref{Orlicz_Hardy}) holds for $f$. 
Therefore, following the terminology of \cite{LQR},  $T$ is uniquely extended as a bounded operator from the Orlicz-Hardy space $H^{\Phi_r} (\mathbb{T})$ to $L^1 (\mathbb{T})$, where  $\Phi_r (x) = x ( [1+\log (x+1)]^r - 1 )$, $x \geq 0$.

\subsection{Higher-dimensional variants}\label{open_problems}
In \cite{BRS}, both the results of Zygmund and Pichorides were extended to higher dimensions. Hence, in analogy to the one-dimensional case, one is naturally led to ask whether versions of Theorem \ref{main_theorem} involving Hardy spaces of several variables hold as well. However, it should be mentioned that it does not seem that the methods of this paper can be extended in a straightforward way to operators acting on Hardy spaces in polydiscs. 

\subsection*{Acknowledgements.} The author would like to thank Salvador Rodr\'iguez-L\'opez and Alan Sola for some interesting discussions and for carefully reading an earlier version of this manuscript.

\end{document}